%

\input amstex
\documentstyle{amsppt}

\UseAMSsymbols

\catcode`\@=11

\def\leftrightarrowfill{$\m@th\mathord\leftarrow\mkern-6mu%
  \cleaders\hbox{$\mkern-2mu\mathord-\mkern-2mu$}\hfill
  \mkern-6mu\mathord\rightarrow$}

\atdef@?#1?#2?{\ampersand@\setbox\z@\hbox{$\ssize
 \;\;{#1}\;$}\setbox\@ne\hbox{$\ssize\;\;{#2}\;$}\setbox\tw@
 \hbox{$#2$}\ifCD@
 \global\bigaw@\minCDaw@\else\global\bigaw@\minaw@\fi
 \ifdim\wd\z@>\bigaw@\global\bigaw@\wd\z@\fi
 \ifdim\wd\@ne>\bigaw@\global\bigaw@\wd\@ne\fi
 \ifCD@\hskip.5em\fi
 \ifdim\wd\tw@>\z@
 \mathrel{\mathop{\hbox to\bigaw@{\leftrightarrowfill}}\limits^{#1}_{#2}}\else
 \mathrel{\mathop{\hbox to\bigaw@{\leftrightarrowfill}}\limits^{#1}}\fi
 \ifCD@\hskip.5em\fi\ampersand@}
\atdef@-#1-#2-{\ampersand@\setbox\z@\hbox{$\ssize  
 \;\;{#1}\;$}\setbox\@ne\hbox{$\ssize\;\;{#2}\;$}\setbox\tw@
 \hbox{$#2$}\ifCD@
 \global\bigaw@\minCDaw@\else\global\bigaw@\minaw@\fi
 \ifdim\wd\z@>\bigaw@\global\bigaw@\wd\z@\fi
 \ifdim\wd\@ne>\bigaw@\global\bigaw@\wd\@ne\fi
 \ifCD@\hskip.5em\fi
 \ifdim\wd\tw@>\z@
 \mathrel{\mathop{\raise.5ex\hbox to\bigaw@{\hrulefill}}\limits^{#1}_{#2}}
\else\mathrel{\mathop{\raise.5ex\hbox to\bigaw@{\hrulefill}}\limits^{#1}}\fi
 \ifCD@\hskip.5em\fi\ampersand@}

\def\hookrightarrowfill{$\m@th\mathord\lhook\mkern-3.1mu%
  \cleaders\hbox{$\mkern-2mu\mathord-\mkern-2mu$}\hfill
  \mkern-6mu\mathord\rightarrow$}
\atdef@(#1(#2({\ampersand@\setbox\z@\hbox{$\ssize
 \;\;{#1}\;$}\setbox\@ne\hbox{$\ssize\;\;{#2}\;$}\setbox\tw@
 \hbox{$#2$}\ifCD@
 \global\bigaw@\minCDaw@\else\global\bigaw@\minaw@\fi
 \ifdim\wd\z@>\bigaw@\global\bigaw@\wd\z@\fi
 \ifdim\wd\@ne>\bigaw@\global\bigaw@\wd\@ne\fi
 \ifCD@\hskip.5em\fi
 \ifdim\wd\tw@>\z@
 \mathrel{\mathop{\hbox to\bigaw@{\hookrightarrowfill}}\limits^{#1}_{#2}}\else
 \mathrel{\mathop{\hbox to\bigaw@{\hookrightarrowfill}}\limits^{#1}}\fi
 \ifCD@\hskip.5em\fi\ampersand@}

\parindent20\p@
\advance\captionwidth@1in

\font\chapheadfont@=cmr10 scaled\magstep2
\font\chapheadmathfont@=cmmi10 scaled\magstep2
\font\chapheadmsa=msam10 scaled\magstep2
\font\chapheadmsb=msbm10 scaled\magstep2
\font\chapheadscriptmathfont=cmmi7 scaled\magstep2
\define\chapheadfont{\chapheadfont@\textfont1=\chapheadmathfont@%
        \scriptfont1=\chapheadscriptmathfont%
        \textfont\msafam=\chapheadmsa\textfont\msbfam=\chapheadmsb}
\outer\def\chapheading{\newpage
  \begingroup\raggedcenter@\interlinepenalty\@M \let\\\linebreak
  \chapheadfont\noindent\ignorespaces}

\newif\ifdatver
\define\datver#1{\ifdatver\global\def\datverp{}\else\datvertrue%
  \global\def\datverp{\par\boxed{\boxed{\text{Version: #1; Run: \today}}}}\fi}

\newif\ifbigdoc

\let\section\relax

\def\nosection{\newcodes@\endlinechar=10 \sect@}
{\lccode`\!=`\\
\lowercase{\gdef\sect@#1^^J{\sect@@#1!section\sect@@@}%
\gdef\sect@@#1!section{\futurelet\next\sect@@@}%
\gdef\sect@@@#1\sect@@@{\ifx\next\sect@@@\let
\next=\sect@\else\def\next{\oldcodes@\endlinechar=`\^^M\relax}%
 \fi\next}}}
\catcode`\@=\active
\define\protag#1 #2{{\hbox{\rm\ignorespaces#1\unskip}}#2}
\define\theprotag#1 #2{#2 {\rm\ignorespaces#1\unskip}}
\define\exertag #1 #2{\demo{\hbox{\rm(\ignorespaces#1\unskip)} #2}}

\define\figtagg#1{Figure $#1$}

\NoRunningHeads
\CenteredTagsOnSplits
\NoBlackBoxes

\def\today{\ifcase\month\or
 January\or February\or March\or April\or May\or June\or
 July\or August\or September\or October\or November\or December\fi
 \ \space\number\day, \number\year}

\define\OS {\negthinspace\smallsetminus\negthinspace 0}

\define\dCc#1{\dot{\Cal{C}}^{#1}_c}

\define\nov#1{{\frac{n}{#1}}}

\define\({\bigl(}
\define\){\bigr)}
\hsize=6.5truein
\vsize=8.5truein

\define\cof{\operatorname{cof}}
\define\Sat{\operatorname{Sat}}
\define\Seq{\operatorname{Seq}}
\define\Lim{\operatorname{Lim}}

\document
\baselineskip=.25truein

\font\bigtenrm=cmr12 scaled\magstep2
\centerline
{\bigtenrm {THE GENERICITY CONJECTURE}}

\font\bigtenrm=cmr10 scaled\magstep2
\centerline
{\bigtenrm {Sy D. Friedman\footnote"*"{Research Supported by the National
Science 
Foundation,  Grant \#8903380--DMS.\hphantom{XXXXXXXXXXXXXXXXXXXXXXXX}} }}

\font\bigtenrm=cmr10 scaled\magstep2
\centerline
{\bigtenrm {MIT}}

\vskip20pt

The Genericity Conjecture, as stated in Beller-Jensen-Welch [82], is the
following:
$$\text{ If } O^\#\notin L[R],\ R\subseteq \omega\text{ then } 
R\text{ is generic over } L.\tag{$*$}$$
We must be precise about what is meant by ``generic''. 

\vskip5pt

\flushpar
{\bf Definition.} \ (Stated in Class Theory) \ A {\it generic extension} of an
inner model $M$ is an inner model $M[G]$  such that for some forcing notion
${\Cal{P}}\subseteq M:$ 

\roster
\item"(a)" \ $\langle M,{\Cal{P}}\rangle$ is amenable and $\Vdash_p$ is
$\langle M,{\Cal{P}}\rangle$-definable for $\underset\sim\to\Delta_0$ 
sentences.

\item"(b)" \ $G\subseteq {\Cal{P}}$ is compatible, closed upwards and
intersects every $\langle M,{\Cal{P}}\rangle$-definable dense $D\subseteq
{\Cal{P}}.$  
\endroster

A set $x$  is {\it generic} over $M$  if it is an element of a  generic
extension of $M.$  And $x$ is {\it strictly generic} over $M$  if $M[x]$  is a
generic extension of $M.$ 

Though the above definition quantifies over classes, in the special case where
$M=L$  and $O^\#$  exists these notions are in fact first-order, as all
$L$-amenable classes are $\underset\sim\to\Delta_1$ 
definable over $L[O^\#].$
>From now on assume that $O^\#$ exists.

\vskip5pt

\proclaim{Theorem A} The Genericity Conjecture is false.
\endproclaim

The proof is based upon the fact that every real generic over $L$  obeys a
certain definability property, expressed as follows.

\vskip5pt

\flushpar
{\bf Fact.} \ If $R$ is generic over $L$ then for some $L$-amenable class $A,$
Sat$\langle L,A\rangle$ is {\it not} definable over $\langle L[R], A\rangle,$
where Sat$\langle L,A\rangle$  is the canonical satisfaction predicate for
$\langle L,A\rangle.$ 

Thus Theorem A is established by producing a real $R$ s\.t\. $O^\#\notin L[R]$
yet Sat$\langle L,A\rangle$ {\it is} definable over $\langle L[R],A\rangle$
for each $L$-amenable $A.$ 

A weaker version of the Genericity Conjecture would state: \ If $O^\#\notin
L[R]$  then either $R\in L$  or $R$  is generic over some inner model $M$  not
containing $R.$  
This version of the conjecture is still open. However, this question can also
be studied in contexts where $O^{\#}$ does {\it not} exist, for example when
the universe has ordinal height equal to that of the minimal transitive model
of $ZF.$  In the latter context, Mack Stanley [93] has demonstrated the
consistency of the existence of a non-constructible real which belongs to
every inner model over which it is generic.

\newpage

\flushpar
{\bf Section A \ A Non-Generic Real below $\bold{O^\#}.$ }

We first prove the Fact stated in the introduction.

\proclaim{Lemma 1} \  Suppose $R\subseteq \omega$ 
is generic over $L.$  Then for
some $L$-amenable class $A,$  Sat$\langle L,A\rangle$ is {\bf not} definable
over $\langle L[R],A\rangle$  with parameters.
\endproclaim

\demo{Proof} Let $R\in L[G]$ where $G\subseteq {\Cal{P}}$ is generic for
$\langle L,{\Cal{P}}\rangle$-definable dense classes and ${\Cal{P}}$ is
$L$-amenable as 
in (a), (b) of the definition of generic extension. Let $A={\Cal{P}}$ and
suppose that Sat$\langle L,{\Cal{P}}\rangle$ were definable over $\langle
L[R],{\Cal{P}}\rangle$ with parameters. But the Truth Lemma holds for
$G,{\Cal{P}}$ for formulas mentioning $G,{\Cal{P}}:\langle
L[G],G,{\Cal{P}}\rangle\vDash\phi(G,{\Cal{P}})$  iff $\exists p\in
G(p\Vdash\phi(G,{\Cal{P}})),$ using the fact that $\Vdash$ in ${\Cal{P}}$ for
$\underset\sim\to\Delta_0$ 
sentences is definable over $\langle L,{\Cal{P}}\rangle$  and the
genericity of $G.$  So Sat$\langle L[G],G,{\Cal{P}}\rangle$ is definable over
$\langle L[G], G,\text{Sat}\langle L,{\Cal{P}}\rangle\rangle,$ since $\Vdash$
is definable over $\langle L,\Sat\langle L,{\Cal{P}}\rangle\rangle$ for
arbitrary first-order sentences. Since Sat$\langle L,{\Cal{P}}\rangle$ is
definable over $\langle L[G], G,{\Cal{P}}\rangle$ we get the definability of
satisfaction
for the latter structure over itself. This contradicts a well-known
result of Tarski. \hfill{$\dashv$}
\enddemo

The rest of this section is devoted to the construction of a real $R$  such
that $R$  preserves $L$-cofinalities (cof$(\alpha)$ in $L=\cof(\alpha )$ in
$L[R]$ for every $\alpha$) and for every $L$-amenable $A,$ $\Sat\langle
L,A\rangle$ is definable over $\langle L[R],A\rangle.$ (The proof has little
to do with the Sat operator; any operator from $L$-amenable classes to
$L$-amenable classes that is ``reasonable'' is codable by a real. We discuss
this further at the end of this section.)

$R$  will generically code a class $f$  which is generic for a forcing of size
$\infty^+=$ least  ``$L$-cardinal'' greater than $\infty.$  Since this sounds
like nonsense we suggest that the reader think of $\infty$ as some uncountable
cardinal of $V$  and then $\infty^+$ denotes $(\infty^+)^L.$ Thus we will
define a constructible set forcing ${\Cal{P}}^\infty\subseteq L_{\infty^{+}}$
for adding a generic $f^\infty\subseteq \infty$ such that if
$A\subseteq\infty$  is constructible then Sat$\langle L_\infty,A\rangle$ is
definable over $\langle L_\infty [f^\infty ],f^\infty,A\rangle.$  Then we show
how to choose the $f^\infty$'s to ``fit together'' into an $f\subseteq ORD$
such that Sat$\langle L,A\rangle$ is definable over $\langle L[f],f,A\rangle$
for each $L$-amenable $A.$  Finally, we code $f$  by a real $R$  (using the
fact that $I=$ Silver Indiscernibles are indiscernibles for $\langle L[f],
f\rangle).$  

A condition in ${\Cal{P}}^\infty$ is defined as follows. Work in $L.$  An {\it
Easton set of ordinals} is a set of ordinals $X$  such that $X\cap\kappa$ is
bounded in $\kappa$ for every regular $\kappa>\omega.$  For any $\alpha\in
ORD,$ $2^\alpha$ denotes all $f:\ \alpha\longrightarrow 2$ and
$2^{<\alpha}=\cup\{2^\beta|\beta<\alpha\}.$  An {\it Easton set of strings} is
a 
set $D\subseteq \cup\{2^\alpha|\alpha\in ORD\}$ such that
$D\cap2^{<\kappa }$ has cardinality less than $\kappa$  for every
regular $\kappa>\omega.$  For any $X\subseteq ORD$ let
Seq$(X)=\cup\{2^\alpha|\alpha\in X\}.$  A {\it condition} in
${\Cal{P}}^\infty$ is $(X,F,D,f)$  where:
\roster
\item"(a)" \ $X\subseteq\infty$ is an Easton set of ordinals

\item"(b)" \ $F:\ X\longrightarrow {\Cal{P}}(2^\infty)=\text{ Power Set of
}2^\infty$ such that for $\alpha\in X,$ $F(\alpha)$ has cardinality
$\le\alpha$ 

\item"(c)" \ $D\subseteq\Seq(X)$ is an Easton set of strings

\item"(d)" \ $f:\ D\longrightarrow\infty$ such that $f(s)>\text{ length
}(s)\text{ for } s\in D.$ 
\endroster
We define extension of conditions as follows. $(Y,G,E,g)\le(X,F,D,f)$ iff

\roster
\item"(i)" \ $Y\supseteq X,\ E\supseteq D,\ G(\alpha)\supseteq F(\alpha)$ for
$\alpha\in X, g$ extends $f$ 

\item"(ii)" \ If $s\in E-D$ then the interval (length $(s)+1, g(s)$] contains no
element of $X,$  and if $s\subseteq S\in F(\alpha)$ for some $\alpha\le$
length $(s),\alpha\in X$ then $g(s)\notin C_S.$ 
\endroster
We must define $C_S.$  For $S\in 2^\infty$ let $\mu(S)=$ least p\.r\. closed
$\mu>\infty$ such that $S\in L_\mu$ and then
$C_S=\{\alpha<\infty|\alpha=\infty\cap$ Skolem hull $(\alpha)$
in $L_{\mu(S)}\}.$  Thus $C_S$ is  CUB in $\infty$ and $\langle
L_\alpha,S\restriction\alpha\rangle\prec\langle L_\infty,S\rangle$ for
sufficiently large $\alpha\in C_S$ (as $S\in$ Skolem hull $(\alpha)$ in
$L_{\mu(S)}$ for sufficiently large $\alpha<\infty).$ 
Also note that $T\notin L_{\mu(S)}\longrightarrow C_T\subseteq \Lim
C_S\cup\alpha$ for some $\alpha<\infty.$  

Our goal with this forcing is to produce a generic function $f_G$ from
$2^{<\infty }$  into $\infty$ such that for each $S\subseteq\infty,$
$\{f(S\restriction\alpha) | \alpha < \infty\}$ is a good approximation to
the complement of $C_S$. $S\in F(\alpha)$ is a committment
that for $\beta>\alpha,f(S\restriction\beta)\notin C_S$ (in stronger
conditions).

\proclaim{Lemma 2} If $p\in {\Cal{P}}^\infty$  and $\alpha<\infty,S\in
2^\infty, s\in 2^{<\infty }$  then $p$  has an extension $(X,F,D,f)$  such
that $\alpha\in X,$  $S\in F(\alpha)$  and $s\in D.$ 
\endproclaim

\demo{Proof} Easy, given the fact that if $s$ needs to be added then we can
safely put $f(s)=\text{ length}(s)+1.$ \hfill{$\dashv$}
\enddemo

\proclaim{Lemma 3} ${\Cal{P}}^\infty$ has the $\infty^+$-chain condition
(antichains have size $\le\infty,$  all in $L$ of course).
\endproclaim

\demo{Proof} Any two conditions $(X,F,D,f),(X,G,D,f)$  are compatible, so an
antichain has cardinality at most the number of $(X,D,f)$'s, which is
$\infty.$  \hfill{$\dashv$}
\enddemo

\proclaim{Lemma 4} Let $G$  be ${\Cal{P}}^\infty$-generic and write $f_G$ for
$\cup\{f|(X,F,D,f)\in G$  for some $X,F,D\}.$  If $S\in 2^\infty$ then
$f_G(S\restriction\alpha)\notin C_S$ for sufficiently large $\alpha<\infty.$  
\endproclaim

\demo{Proof} $G$ contains a condition $(X,F,D,f)$  such that $0\in X$  and
$S\in F(0).$  If $s\subseteq S, s\notin D$ then $f_G(s)\notin C_S,$  by (ii)
in the definition of extension. And $S\restriction\alpha\notin D$  for
sufficiently large $\alpha<\infty.$  \hfill{$\dashv$}
\enddemo

\proclaim{Lemma 5} Let $G,f_G$ 
be as in Lemma 4. If $\alpha<\infty$ is regular,
$S\in 2^\infty,$ and $\alpha\notin Lim$ $C_S$ then
$\{f_G(S|\beta)|\beta < \alpha\}$ intersects every constructible
unbounded subset of $\alpha.$
\endproclaim

\demo{Proof} Let $A\subseteq \alpha$ be constructible and unbounded in
$\alpha$. We show that a condition $(X,F,D,f)$ can be extended to
$(X\cup\{\delta\},F^*,D\cup\{S\restriction\delta\},f^*)$ 
for some $\delta,$ where
$f^*(S\restriction\delta)\in A$.
Choose $\delta<\alpha$ large enough so that
$S\restriction\delta$ is not an initial segment of any
$T\in \cup\{F(\beta)|\beta\in X\cap\alpha\}-\{S\}.$ This is possible since
$X\cap\alpha$ is bounded in $\alpha$  and $F(\beta)$ has cardinality $<\alpha$
for each $\beta\in X\cap\alpha.$  Then let $f^*=f\cup\{\langle
S\restriction\delta,\beta\rangle\}$ where $\beta\in A - C_S - \delta$
and $F^*=F\cup\{\langle\delta,\emptyset\rangle\}.$  \hfill{$\dashv$}
\enddemo

\proclaim{Lemma 6} ${\Cal{P}}^\infty$ preserves cofinalities (i\.e\.,
${\Cal{P}}^\infty\Vdash\cof(\alpha)=\cof(\alpha)$ in $L$ for every ordinal
$\alpha$). 
\endproclaim

\demo{Proof} For regular $\kappa<\infty$ and $p\in{\Cal{P}}^\infty$ let
$(p)^\kappa =$ ``part of $p$  below $\kappa$'', $(p)_\kappa=$ ``part of $p$
at or above $\kappa$'' be defined in the natural way: if  $p=(X,F,D,f)$
then
$$(p)^\kappa=(X\cap\kappa,F\restriction
X\cap\kappa,D\cap\Seq\kappa,f\restriction D\cap\Seq\kappa)\text{ and }$$
$$(p)_\kappa=(X-\kappa,F\restriction X-\kappa,D\cap\Seq(\infty-\kappa),f
\restriction D\cap\Seq(\infty-\kappa)).$$
Given $p$  and predense $\langle\Delta_i|i<\kappa\rangle$ we find $q\le p$ and
$\langle\overline\Delta_i|i<\kappa\rangle$ such that
$\overline\Delta_i\subseteq\Delta_i$  for all $i<\kappa,$
card $\overline\Delta_i\le\kappa$  for all $i<\kappa$  and each
$\overline\Delta_i$ is predense below $q.$  ($\Delta$ is {\it predense} if
$\{r|r\le$ some $d\in\Delta\}$ is dense; it is {\it predense} {\it below} $q$
if every extension  of $q$  can be extended into the afore-mentioned set.)
This implies that if $\cof(\alpha)\le \kappa$ in some generic extension
$L[G],G \ {\Cal{P}}^\infty$-generic over $L,$  then $\cof(\alpha)\le\kappa$ in
$L.$  Since ${\Cal{P}}^\infty$ is $\infty^+$-CC, this means that
${\Cal{P}}^\infty$  preserves all cofinalities.

Given $p$  and $\langle\Delta_i|i<\kappa\rangle$  as above first extend $p$ to
$p_0=(X_0, F_0, D_0, f_0)$ so that $\kappa\in X_0.$  Now note that if $r\le
p_0$  then $f^r(s)<\kappa$  for all $s\in D^r-D_0$  of length $<\kappa$ (where
$r=(X^r,F^r,D^r,f^r)),$  by condition (ii) in the definition of extension. 
Thus ${\Cal{F}}=\{(X^r\cap\kappa,D^r\cap\Seq\kappa,f^r\restriction
D^r\cap\Seq\kappa)|r\le p_0\}$ is a set of cardinality $\kappa.$  Let
$\langle(\Delta^*_i,(X^i, D^i, f^i))|i<\kappa)$  be an enumeration in length
$\kappa$  of all pairs from $\{\Delta_i|i<\kappa\}\times {\Cal{F}}.$ 

Now we extend $p_0$  successively to $p_1\ge p_2\ge\dots$ in $\kappa$ steps so
that $(p_i)^\kappa=(p_0)^\kappa$  for all $i<\kappa,$  according to the
following prescription: \ If $p_i$  has been defined, see if it has an
extension $r_i$  extending some $d_i\in\Delta^*_i$  such that $(X^{r_{i}}
\cap\kappa, D^{r_{i}}\cap\Seq\kappa,f^{r_{i}}\restriction
D^{r_{i}}\cap\Seq\kappa)=(X^i,D^i, f^i).$  If not then $p_{i+1}=p_i.$  If so,
select such an $r_i, d_i$  and define $p_{i+1}$  by requiring
$(p_{i+1})^\kappa=(p_0)^\kappa,(p_{i+1})_\kappa=(r_i)_\kappa$ except enlarge
$F^{p_{i+1}}(\kappa)$ so as to contain  $F^{r_{i}}(\alpha)$ for $\alpha\in
X^{r_{i}}\cap\kappa.$  For limit $\lambda\le\kappa$  let $p_\lambda$  be
the greatest lower bound to $\langle p_i|i<\lambda\rangle.$  Finally let
$q=p_\kappa.$  

Let $\overline\Delta_j\subseteq\Delta_j$  consist of all $d_i$ in the above
construction that belong to $\Delta_j,$  for $j<\kappa.$  The claim we must
establish is that each $\overline\Delta_j$ is predense below $q.$  Here's the
proof: suppose $\bar q\le q$ and  let $r\le\bar q, r$ extending some element
of $\Delta_j.$ Choose $i<\kappa$ so that
$(\Delta^*_i,(X^i,D^i,f^i))=(\Delta_j,(X^r\cap\kappa,D^r\cap\Seq\kappa,f^r
\restriction D^r\cap\Seq\kappa)).$   Clearly at stage $i+1,$ it was possible
to find $r_i, d_i$  as searched for in the construction. It suffices to argue
that $r_i,\bar q$  are compatible. Now $(r_i)_\kappa$  is extended by
$(p_{i+1})_\kappa$  and hence by $(r)_\kappa.$  And $(r_i)^\kappa$ is extended
by $(r)^\kappa,$  except possibly that $F^{r_{i}}(\alpha)$  may fail to be a
subset of $F^r(\alpha)$  for $\alpha\in X^r\cap\kappa.$  And note that the
extension $(r_i)_\kappa\ge(r)_\kappa$  obeys all restraint imposed by
$F^{r_{i}}(\alpha)$  for $\alpha\in X^r\cap\kappa$ since we included
$F^{r_{i}}(\alpha)$  in $F^{p_{i+1}}(\kappa).$  Thus $r_i$  and $\bar q$  are
both extended by $r,$  provided we only enlarge $F^r(\alpha)$ for $\alpha\in
X^r\cap\kappa$ to include $F^{r_{i}}(\alpha),$ \hfill{$\dashv$}
\enddemo

For future reference we state:

\proclaim{Corollary 6.1} Suppose $\kappa<\infty$ is regular and
$\Delta\subseteq {\Cal{P}}^\infty$ is predense. 
Let ${\Cal{P}}^\infty_\kappa=\{(p)_\kappa|p\in{\Cal{P}}^\infty\},
{\Cal{P}}^{\infty,\kappa }=\{p\in{\Cal{P}}^\infty|X^p\subseteq\kappa$ and
Range $(f^p)\subseteq\kappa\}$  with the notion $\le$  of extension defined as
for ${\Cal{P}}^\infty.$  Then for any $q\in {\Cal{P}}^\infty_\kappa$  there is
$q'\le q$ such that $\Delta^{q'}=\{r\in {\Cal{P}}^{\infty,\kappa }|r\cup q'$
meets $\Delta,F^r(\alpha)\subseteq F^{q'}(\kappa)$ 
for all $\alpha\in X^r\}$ is predense on ${\Cal{P}}^{\infty,\kappa}.$ 
\endproclaim

\demo{Proof} As in the proof of Lemma 6, successively extend $q$  (after
guaranteeing $\kappa\in X^q)$ in $\kappa$  steps to $q'$  so that for any
$(X,D,f)$ if  $r\cup q''$ meets $\Delta$  for some $q''\le q',$  some $r$
such that $(X^r,D^r,f^r)=(X,D,f)$  then $r\cup q'$  meets $\Delta$  for some
such $r,$  where $F^r(\alpha)\subseteq F^{q'}(\kappa)$  for all $\alpha\in
X^r.$  Now note that if $r_0\in{\Cal{P}}^{\infty,\kappa }$  then $r_0\cup q'$
has 
an extension meeting $\Delta$  so there is $r_1$ such that
$(X^{r_{1}},D^{r_{1}}, f^{r_{1}})=(X^{r_{0}},D^{r_{0}}, f^{r_{0}})$  and
$r_1\in\Delta^{q'}.$  But then $r_0$  is compatible with $r_1$  so
$\Delta^{q'}$  is predense on ${\Cal{P}}^{\infty,\kappa },$  as desired.
\hfill{$\dashv$} 
\enddemo

\proclaim{Corollary 6.2} ${\Cal{P}}^\infty\Vdash GCH.$ 
\endproclaim

\demo{Proof} Suppose $f^\infty:\ \Seq(\infty)\longrightarrow \infty$ is
${\Cal{P}}^\infty$-generic. It suffices to show that if $\kappa\le\infty$ is
regular, $A\subseteq\kappa,$ $A\in L[f^\infty ]$ then $A\in
L[f^\infty\restriction\Seq(\kappa)].$  But the proof of Lemma 6 shows that
given any $p\Vdash\dot A\subseteq\kappa$ there is $q\le p$ such that for any
$i<\kappa,$  $\{r\le q|(r)_\kappa=(q)_\kappa$  and $r$  decides 
``$i\in\dot A\text{''}\}$ is predense below $q.$  This proves that there is
$q\le p$  such that $q\Vdash\dot A\in L[\dot
f^\infty\restriction\Seq(\kappa)]$  and so by the genericity of $f^\infty,$
$A\in L[f^\infty\restriction\Seq(\kappa)].$  \hfill{$\dashv$}
\enddemo

Next we embark on a series of lemmas aimed at showing that
${\Cal{P}}^\infty$-generics actually exist in $L[O^\#]$ when $\infty$ is any
Silver indiscernible.

\proclaim{Lemma 7} Suppose $i<j$ are adjacent countable Silver indiscernibles.
Let $\pi=\pi_{ij}$ denote the elementary embedding $L\longrightarrow L$ which
shifts each of the Silver indiscernibles $\ge i$  to the next one and leaves
all other Silver indiscernibles fixed. Then there is a ${\Cal{P}}^j_i$-generic
$G^j_i$  such that if $(X,F,D,f)$  belongs to $G^j_i$  and $S\subseteq i,$
$S\in L$  then $f(\pi(S)\restriction\alpha)\notin C_{\pi(S)}$ for all
$\pi(S)\restriction\alpha\in D.$ 
\endproclaim

\demo{Proof} For any $k\in \omega$ let $\ell_1<\dots<\ell_k$ be the first $k$
Silver indiscernibles greater than $j$  and let $j_k=j^+\cap\Sigma_1$ Skolem
hull of $j+1\cup\{\ell_1\dots\ell_k\}$ in $L,$ $i_k=i^+\cap\Sigma_1$ Skolem
hull of $i+1\cup\{\ell_1\dots\ell_k\}$ in $L.$  (Of course $i^+,j^+$ denote
the cardinal successors to $i,j$ {\it in} $L$.) Let $j^*_k=$ least p\.r\.
closed ordinal $\alpha>j_k$ such that $L_\alpha\vDash j$ is the largest
cardinal. Finally let $C_k=\{\gamma<j|\gamma=j\cap\Sigma_1$ Skolem hull
$(\gamma\cup\{j\}\cup\{\ell_1\dots\ell_k))$ in $L\},$ a CUB subset of $j.$ 

Now note that if $S\subseteq i, S\in L-L_{i_{k}}$ then $C_{\pi(S)}\subseteq
C_k\cup\gamma$  for some $\gamma<i.$  For, $\mu_{\pi(S)}$ is greater than or
equal to $j^*_k$ since otherwise $\pi(S)$ belongs to $L_{j_{k}}$  and hence
$S$  belongs to $L_{i_{k}}.$  Thus $C_{\pi(S)}\subseteq C_k\cup\gamma$ for
some $\gamma<j$  since $C_k$  is an element of $L_{j^{*}_{k}}.$ But the least
such $\gamma$ is definable from elements of $i\cup$ (Silver Indiscernibles
$\ge j),$  so must be less than $i.$ 

Also note that the $L$-cofinality of $j_k$ is equal to $j:$  Consider
$M=$transitive collapse of $\Sigma_1$  Skolem hull of
$j+1\cup\{\ell_1\dots\ell_k\}.$  There is a partial $\underline{\Sigma_1}(M)$ 
function from a subset of $j$  onto $j_k,$  all of whose restrictions to
ordinals $\gamma<j$  have range bounded in $j_k.$  (This is why we are using
$\Sigma_1$  Skolem hulls rather than full $\Sigma_\omega$ Skolem hulls.) Thus
the $L$-cofinalities of $j_k$  and $j$  are the same, namely $j.$ 

Thus we may conclude the following: The set $\{\pi(S)|S\subseteq i, S\in
L_{i_{k}}\}\in L_{j_{k}}$ (since it is a constructible bounded subset of
$L_{j_{k}}$) 
and if $S\subseteq i,$ $S\in L-L_{i_{k}}$  then
$C_{\pi(S)}$ is disjoint from $(i,\gamma _k),$ where $\gamma_k=$ least
element of $C_k$  greater than $i.$ 

Now we see how to build $G^j_i.$  We describe an $\omega$-sequence $p_0\ge
p_1\ge\dots$ of conditions in ${\Cal{P}}^j_i$  and take $G^j_i=\{p\in
{\Cal{P}}_i^j|p_k\le p$  for some $k\}.$  Let
$\langle\Delta_k|k\in\omega\rangle$  be a list of all constructible dense sets
on ${\Cal{P}}^j_i$ so that for all $k,$  $\Delta_k$ belongs to the $\Sigma_1$
Skolem hull in $L$  of $i\cup\{i,j,\ell_1\dots\ell_{k+1}\}.$  This is possible
since any constructible dense set on ${\Cal{P}}^j_i$  belongs to $L_{j^{++}}$
and hence to the $\Sigma_1$ Skolem hull in $L$  of
$i\cup\{i,j,\ell_1\dots\ell_k\}$  for some $k.$  We inductively define $p_0\ge
p_1\ge\dots$  so that $p_k$  belongs to the $\Sigma_1$  Skolem hull in $L$  of
$i^+\cup\{j,\ell_1\dots\ell_k\}.$  Let $p_0$  be the weakest condition in
${\Cal{P}}^j_i; \ p_0=(\emptyset,\emptyset,\emptyset,\emptyset).$ Suppose that
$k>0$ and $p_{k-1}$  has been defined. Write $p_{k-1}=(X,F,D,f).$  First
obtain $\bar p_k$  by adding $i$  to $X$ if necessary and defining or
enlarging $F(i)$  so as to include  
$\{\pi(S)|S\subseteq i, S\in L_{i_{k}}\}.$  Then choose $p_k\le\bar p_k$ to be
$L$-least so that $p_k$  meets $\Delta_{k-1}.$ This completes the
construction. 

We show that $p_k\in\Sigma_1$  Skolem hull in $L$  of
$i^+\cup\{j,\ell_1\dots\ell_k\}.$  By induction $p_{k-1}$  belongs to this
hull and by choice of $\langle\Delta_k|k\in\omega\rangle$, so does
$\Delta_{k-1}.$ Now $\{\pi(S)|S\subseteq i, S\in L_{i_{k}}\}$ is the range of
$f\restriction i$  where $f$  is a $\underline{\Sigma_1}(L)$ partial function
with parameters $j,\ell_1\dots\ell_k.$  The latter is because
Range$(\pi\restriction i_k)$ is just $j_k\cap\Sigma_1$ Skolem hull in $L$
of $i\cup\{j,\ell_1\dots\ell_k\}.$  But given a parameter $x$  for the domain
of this $\underline{\Sigma_1}(L)$ partial function, its range becomes 
$\Sigma _1$-definable in the sense that it is in the $\Sigma_1$ Skolem hull in
$L$ of $\{x,j,\ell_1\dots\ell_k\}.$  As $x$  can be chosen equivalently as an
ordinal $<i^+,$  we get that $\{\pi(S)|S\subseteq i,S\in L_{i_{k}}\}$ belongs
to the $\Sigma_1$  Skolem hull in $L$ of $i^+\cup\{j,\ell_1\dots\ell_k\}.$
Thus so does $p_k.$  (Actually $x$  can be chosen to be $i_k.$)

Finally we must check that if $p_k=(X_k,F_k,D_k,f_k)$ then
$f_k(\pi(S)\restriction\alpha)\notin C_{\pi(S)}$ for all
$\pi(S)\restriction\alpha\in D_k,$  all $S\subseteq i$ in $L.$  Assume that
this is true for smaller $k$  and we check it for $k.$  Now if $S\in
L_{i_{k}}$ then this is guaranteed by the fact that
$\pi(S)\in\overline{F}_k(i),$  where $\bar
p_k=(\overline{X}_k,\overline{F}_k,D_{k-1},f_{k-1}).$  If $S\in L-L_{i_{k}}$
then $C_{\pi(S)}$ is disjoint from $(i,\gamma_k)$, where
$\gamma_k=j\cap\Sigma_1$ Skolem hull in $L$  of
$\gamma_k\cup\{j\}\cup\{\ell_1\dots\ell_k\}$  and $\gamma_{k}>i.$  But then
$\gamma_k>i^+$ so $C_{\pi(S)}$ is disjoint from $(i,\bar\gamma_k)$ where
$\bar\gamma_k=\sup(j\cap\Sigma_1$ Skolem hull in $L$  of
$i^+\cup\{j\}\cup\{\ell_1\dots\ell_k\}).$  Since $p_k\in\Sigma_1$ Skolem hull
in $L$ of $i^+\cup\{j\}\cup\{\ell_1\dots\ell_k\},$ it follows that
Range$(f_k)\subseteq\bar\gamma_k$  and hence Range$(f_k)$ is disjoint from
$C_{\pi(S)}.$  \hfill{$\dashv$}
\enddemo
\proclaim{Lemma 8} Suppose $i<j$  are adjacent Silver indiscernibles, $G^j_i$
is ${\Cal{P}}^j_i$-generic over $L$  as in Lemma 7 and $G^i$ is
${\Cal{P}}^i$-generic over $L.$  Then there exists $G^j$  which is
${\Cal{P}}^j$-generic over $L$  such that $G^j_i=\{(p)_i|p\in G^j\}$  and
$q\in G^i\longleftrightarrow\pi_{ij}(q)\in G^j.$ 
\endproclaim

\demo{Proof} As before, let ${\Cal{P}}^{j,i}\subseteq {\Cal{P}}^j$ 
consist of all
$p=(X^p,F^p,D^p,f^p)$ in ${\Cal{P}}^j$  such that $X^p\subseteq i$  and Range
$(f^p)\subseteq i.$  For any $p\in {\Cal{P}}^{j,i}$  we modify $p$  to $\bar
p$  as follows. For $S\in F^p(\alpha),$  $i\in C_S$  let $\bar
S=\pi_{ij}(S\restriction i).$  For $S\in F^p(\alpha), i\notin C_S$
let $T\subseteq i$ be $L$-least so that  $(T,C_T),$ $(S,C_S)$ agree through
sup$(C_S\cap i)$  and let
$\overline S=\pi_{ij}(T).$  Then $F^{\bar p}(\alpha)$  consists of all
$\overline{S}$  for $S\in F^p(\alpha).$  Otherwise $p,\bar p$  agree: $(X^p,
D^p, f^p)=(X^{\bar p}, D^{\bar p}, f^{\bar p}).$ 

If $p\in{\Cal{P}}^j_i$ and $i\in X^p$  we let $Q(p)$  denote
$\{q\in{\Cal{P}}^{j,i}|F^q(\alpha)\subseteq F^p(i)$  for all $\alpha\in
X^q.\}$ Now define $\overline{G}^j=\{p\in{\Cal{P}}^j|(p)_i\in G^j_i, i\in
X^p,$ $(p)^i\in Q((p)_i)$  and $\overline{(p)^i}\in\pi_{ij}[G^i]\}.$  Note
that if $p_0, p_1$  belong to $\overline{G}^j$ then $p_0, p_1$  are
compatible because  $(p_0)_i,(p_1)_i$ are compatible, the restraints from
$(p_0)^i,(p_1)^i$ are ``covered'' by $F^{p_{0}}(i),F^{p_{1}}(i)$ and
$\overline{(p_0)^i},\overline{(p_1)^i}$ 
impose at least as much restraint below
$i$ as do $(p_0)^i,(p_1)^i.$  Note that if $G^j=\{p|\bar p\le p$ for some
$\bar p\in\bar G^j\}$ then $G^j$ is compatible, closed upwards and
$G^j_i=\{(p)_i|p\in G^j\}.$  Also $q\in G^i\longleftrightarrow\pi_{ij}(q)\in
G^j,$  using the hypothesis that $G^j_i$  satisfies Lemma 7. So it only
remains to show that $\overline{G}^j$ meets all constructible predense
$\Delta\subseteq {\Cal{P}}^j.$ 

The first Corollary to Lemma 6 states that it is enough to show that
$\overline{G}^j_i=\{(p)_i|p\in\overline{G}^j\}$ meets all constructible
predense $\Delta\subseteq {\Cal{P}}^j_i$  and that for
$p\in\overline{G}^j_i,\{q\in Q(p)|q=(r)^i$  for some  $r\in\overline{G}^j\}$
meets all constructible $\Delta\subseteq Q(p)$ which are predense on
$\cup\{Q(p^*)|p^*\le p\}={\Cal{P}}^{j,i}.$  
The former assertion is clear by the
${\Cal{P}}^j_i$-genericity over $L$  of $G^j_i=\overline{G}^j_i.$ 
To prove the latter assertion we must show that for $p\in\overline{G}^j_i,$
$\{q\in Q(p)|\overline q\in\pi_{ij}[G^i]\}$  meets every constructible
$\Delta\subseteq Q(p)$ which is predense on ${\Cal{P}}^{j,i}.$  Given such a
$\Delta,$  let $\overline\Delta\subseteq {\Cal{P}}^i$  be defined by
$\overline\Delta=\{r\in {\Cal{P}}^i|\pi_{ij}(r)=\bar q$ for some
$q$ meeting $\Delta\}.$ 
Note that $\overline\Delta$ is constructible because it equals $\{r\in 
{\Cal{P}}^i|r=\pi^{-1}_{ij}(\bar q)$ for some $q$ meeting $\Delta\}$ and
 $\Delta$
has $L$-cardinality $\le i.$  We claim that $\overline\Delta\subseteq
{\Cal{P}}^i$ is predense on ${\Cal{P}}^i.$  Indeed, if $r\in {\Cal{P}}^i$ then
$\pi_{ij}(r)\in {\Cal{P}}^{j,i}$ and therefore can be extended to some $q$ meeting $\Delta.$ As $\bar q=\pi_{ij}(t)$  for some $t\le r$  we have shown
that $r$  can be extended into $\overline\Delta.$  By the
${\Cal{P}}^i$-genericity of $G^i,$  choose $r\in\overline\Delta\cap G^i.$
Then $\pi_{ij}(r)=\bar q$ where $q$ meets $\Delta$; clearly $\bar
q\in\pi_{ij}[G^i].$  \hfill{$\dashv$}
\enddemo

\proclaim{Lemma 9} Let $i_1<i_2<\dots$ denote the first $\omega$-many Silver
indiscernibles and $i_\omega$  their supremum. Then there exist $\langle
G^{i_{n}}|n\ge 1\rangle$ such that $G^{i_{n}}$ is ${\Cal{P}}^{i_{n}}$-generic
over $L$  and whenever $\pi:\ L\longrightarrow L$ is elementary,
$\pi(i_\omega)=i_\omega$ we have $p\in G^{i_{n}}\longleftrightarrow\pi(p)\in
G^{\pi(i_{n})}.$ 
\endproclaim

\demo{Proof} Note that any $\pi$ as in the statement of the lemma restricts to
an increasing map from $\{i_n|n\ge 1\}$ to itself, so $G^{\pi(i_n)}$ makes
sense. We define $G^{i_{n}}$ by induction on $n\ge 1.$  Select $G^{i_1}$ 
to be
the $L[O^\#]$-least ${\Cal{P}}^{i_{1}}$-generic (over $L$). Select
$G^{i_{2}}_{i_{1}}$  as in Lemma 7 and use Lemma 8 to define $G^{i_{2}}$  from
$G^{i_{2}}_{i_{1}},G^{i_{1}}.$  Now suppose that $G^{i_{n}}$  has defined,
$n\ge 2.$  Then define $G^{i_{n+1}}_{i_{n}}$ to be $\{p\in
{\Cal{P}}^{i_{n+1}}_{i_{n}}|\pi_{i_{1}i_{n}}(q)\le p$ for some $q\in
G^{i_{2}}_{i_{1}}\}$  where $\pi_{i_{1}i_{n}}(i_m)=i_{m+n-1}$ for
$m<\omega,\pi_{i_{1}i_{n}}(j)=j$ for $j\in I-i_\omega.$  
Then $G^{i_{n+1}}_{i_{n}}$ is
${\Cal{P}}^{i_{n+1}}_{i_{n}}$-generic, using the $\le i_1$-closure of
${\Cal{P}}^{i_{2}}_{i_{1}}$  and the fact that the collection of constructible
dense subsets of 
${\Cal{P}}^{i_{n+1}}_{i_{n}}$ is the countable union of sets of the form
$\pi_{i_{1}i_{n}}(A),A$ of $L$-cardinality $i_1.$  Moreover
$G^{i_{n+1}}_{i_{n}}$  obeys the condition of Lemma 7 since
$G^{i_{2}}_{i_{1}}$  does and $\pi_{i_{1}i_{n}}$ is elementary. Now define
$G^{i_{n+1}}$  from $G^{i_{n+1}}_{i_{n}},$ $G^{i_{n}}$ using Lemma 8.

To verify $p\in G^{i_{n}}\longleftrightarrow\pi(p)\in G^{\pi(i_n)},$ note that
this depends only on $\pi\restriction L_{i_{\ell}}$  for some $\ell<\omega$
and
any such map is the finite composition of maps of the form $\pi_m,$  where
$\pi_m(i_n)=i_{n+1}$  for $n\ge m,$  $\pi_m(i_n)=i_n$ for $1\le n<m.$  So we
need only verify that for each $m,n,p\in
G^{i_{n}}\longleftrightarrow\pi_m(p)\in G^{\pi_{m}(i_n)}.$ This is trivial
unless $m\le n$  as $m>n\longrightarrow\pi_m(p)=p$ for $p\in
G^{i_{n}}=G^{\pi_{m}(i_{n})}.$ Finally we prove the statement by induction on
$n\ge m.$  If $n=m$ then it follows from the fact that $G^{i_{n+1}}$  was
defined from $G^{i_{n+1}}_{i_{n}},$  $G^{i_{n}}$  so as to obey the conclusion
of Lemma 8. Suppose it holds for $n\ge m$  and we wish to demonstrate the
property for $n+1.$  But $G^{i_{n+1}}$ is defined from
$G^{i_{n+1}}_{i_{n}},G^{i_{n}}$ as $G^{i_{n+2}}$ is defined from
$G^{i_{n+2}}_{i_{n+1}},G^{i_{n+1}}.$  Clearly
$\pi_m[G^{i_{n+1}}_{i_{n}}]\subseteq G^{i_{n+2}}_{i_{n+1}}$  and by induction
$\pi_m[G^{i_{n}}]\subseteq G^{i_{n+1}}.$ Thus $p\in
G^{i_{n+1}}\longrightarrow\pi_m(p)\in G^{\pi_{m}(i_{n+1})}.$ Conversely,
$p\notin G^{i_{n+1}}\longrightarrow p$ incompatible with some $q\in
G^{i_{n+1}}\longrightarrow\pi_m(p)$ incompatible with some $\pi_m(q)\in
G^{\pi_{m}(i_{n+1})}\longrightarrow\pi_m(p)\notin G^{\pi_{m}(i_{n+1})}.$ 
\hfill{$\dashv$}
\enddemo

\proclaim{Lemma 10} There exist $\langle G^i|i\in I\rangle$ such that $G^i$
is ${\Cal{P}}^i$-generic over $L$  and whenever $\pi:\ L\longrightarrow L$ is
elementary, $p\in G^i\longleftrightarrow\pi(p)\in G^{\pi(i)}.$ 
\endproclaim

\demo{Proof} Let $t$ denote a Skolem term for $L;$  thus $L=\{t(j_1\dots
j_n)|t$ a Skolem term, $t$ $n$-ary, $j_1<\dots< j_n$ in $I\}.$ Now define
$t(j_1\dots j_n)\in G^i$ iff $t(\sigma(j_1)\dots\sigma(j_n))\in G^{\sigma(i)}$
where $\sigma$  is the unique order-preserving map from $\{i,j_1\dots j_n\}$
onto an initial segment of $I.$  $(G^i$ for $i<i_w$ is defined in Lemma 9.) We
verify that this is well-defined: \ if $t_1(j_1\dots j_n)=t_2(k_1\dots k_m)$
then let $\sigma^*$ be  the  unique order-preserving map from $\{i,j_1\dots
j_n,k_1\dots k_m\}$ onto an initial segment of $I.$  Then
$t_1(\sigma^*(j_1)\dots\sigma^*(j_n))=t_2(\sigma^*(k_1)\dots\sigma^*(k_m)).$
But $t_1(\sigma^*(j_1)\dots\sigma^*(j_n))\in G^{\sigma^{*}(i)}$ iff
$t_1(\sigma_1(j_1)\dots\sigma_1(j_n))\in G^{\sigma _{1}(i)}$ where $\sigma_1$
is the unique order-preserving map from $\{i,j_1\dots j_n\}$ onto an initial
segment of $I,$ using Lemma 9. The analogous statement holds for $t_2,$ so our
definition is well-defined. The property  $p\in
G^i\longleftrightarrow\pi(p)\in G^{\pi(i)}$ is clear, using our definition.
\hfill{$\dashv$}
\enddemo

Now we are almost done. For any $i\in I$  let $f^i= \cup\{f^p|p\in G^i\}.$  Thus
$f^i:\ 2^{<i}\longrightarrow i.$  And let $f=\cup\{f^i|i\in I\},$  so $f:\
2^{<\infty }\longrightarrow\infty$ ($\infty$ now denotes {\it real} $\infty,$
that is, $\infty=$ ORD).

\proclaim{Lemma 11} (a) \  For any $L$-amenable $A\subseteq$ ORD, SAT$\langle
L,A\rangle$ is definable over $\langle L[f],f,A\rangle.$ 

(b) \ $I$ is a class of indiscernibles for $\langle L[f],f\rangle.$ 

(c)  \ $L[f]\vDash$ GCH.
\endproclaim

\demo{Proof} (a) \ We treat $A$  as an $L$-amenable function $A:\
\infty\longrightarrow 2.$  By Lemmas 4,5 we have that for sufficiently large
$L$-regular $\alpha,\alpha\in Lim$ $C_A\longleftrightarrow$ Range of
$f\restriction\{A\restriction\beta|\beta<\infty\}$ intersects every 
constructible unbounded subset of $\alpha$ (where $C_A$ is defined for
$A$  to be the limit of $C_{A\restriction i},i\in I).$  But for $\alpha$
sufficiently large in $C_A,$ $\langle L_\alpha,A\restriction\alpha\rangle\prec
\langle L,A\rangle$ so Sat$\langle L,A\rangle$ is definable over $\langle
L[f],f,A\rangle.$ 

(b) \ Clear by Lemma 10.

(c) \ By Corollary 6.2. \hfill{$\dashv$}
\enddemo

Finally, using the technique of the proof of Theorem 0.2 of
Beller-Jensen-Welch [82], there is a real $R$  such that $f$ is definable over
$L[R]$  and $I^R=I.$  Thus we conclude.

\proclaim{Theorem 12} There is a real $R\in L[O^\#]$ such that:
\roster
\item"(a)" \ $L,L[R]$ have the same cofinalities

\item"(b)" \ $I^R=I$ 

\item"(c)" \ If $A$  is an $L$-amenable class then Sat$\langle L,A\rangle$ is
definable over $\langle L[R],A\rangle.$  
\endroster

By Lemma 1 we conclude:
\endproclaim

\proclaim{Theorem A} The Genericity Conjecture is false.
\endproclaim

We close this section by mentioning a generalization of the above treatment of
the SAT operator to other operators on classes. For simplicity we first state
our result in terms of $\omega_1,$  rather than $\infty.$ 

\proclaim{Theorem 13} Assume that $O^\#$ exists. Suppose $F$  is a
constructible function from ${\Cal{P}}_L(\omega_1)$ to itself, where
${\Cal{P}}_L(\omega_1)=$ all constructible subsets of (true) $\omega_1.$  Then
there exists a real $R<_LO^\#$ such that $F(A)$ is definable  over $\langle
L_{\omega_{1}}[R],A\rangle$ for all $A\in {\Cal{P}}_L(\omega_1).$ 
\endproclaim

\demo{Proof} Choose $\alpha<\omega_1$ so that $F$ is definable in $L$  from
parameters in $\alpha\cup(I-\omega_1).$  Also we may construct $F',$  defined
from the same parameters, so that for any $A\in {\Cal{P}}_L(\omega_1),$
$F(A)$ is  definable over $\langle L_{\omega_{1}},A,B\rangle$ for any
unbounded $B\subseteq F'(A).$  Finally note that we may assume that
$F'(A)\subseteq C_A$ for all $A$  (where $A$ is viewed as an element of
$2^{\omega _{1}})$  since $C_A$ is definable over $\langle
L_{\omega_{1}},A,B\rangle$ for any unbounded $B\subseteq C_A.$ 

For any $i\in I,$ $\alpha\le i\le\omega_1,$ let $F'_i$ be defined in $L$ just
like $F',$ but with $\omega_1$ replaced by $i.$  Also define ${\Cal{P}}^i$ as
before but with $C_S$ replaced by $F'_i(S)$ (viewing $S\in 2^i$ as a subset of
$i).$  Then as before we can construct a generic
$f:2^{<\omega_1}\longrightarrow \omega _1$ so that for
any $A\in{\Cal{P}}_L(\omega_1),$ $F(A)$ is definable over $\langle
L_{\omega_{1}}[f],A\rangle.$  Finally code $f$  generically by a real using
the fact that $\alpha$  is countable and $I\cap(\alpha,\omega_1)$ is a set of
indiscernibles for $\langle L_{\omega_{1}}[f],f\rangle.$ \hfill{$\dashv$}
\enddemo

To deal with operators on $L$-amenable classes, we have to keep track of
parameters. 

\vskip5pt

\flushpar
{\bf Definition.} \ Suppose $i<j$  belong to $I$  and $F_i$  is a
counstructible function from ${\Cal{P}}_L(i)$ to itself. Then
$F^j_i:{\Cal{P}}_L(j)\longrightarrow {\Cal{P}}_L(j)$ is defined as follows: \
Write $F_i=t(\alpha,i,\vec k)$ where $t$  is a Skolem term for $L,$ $\alpha<i$
and $\vec k$  are Silver indiscernibles greater than $j.$  Then
$F^j_i=t(\alpha,j,\vec k).$ 

Also define $F^\infty_i:\ L$-amenable classes
$={\Cal{P}}_L(\infty)\longrightarrow {\Cal{P}}_L(\infty)$ as follows: \ Given
an $L$-amenable $A$  choose $t$  and $\alpha$ so that for all $j\in I$
greater than $\alpha,$  $A\cap j=t(\alpha,j,\vec k)$ where $\vec k$ are Silver
indiscernibles greater than $j.$  Then $F^\infty_i(A)=\cup\{F^j_i(A\cap
j)|\alpha<j\in I\}.$  An operator $F:\ {\Cal{P}}_L(\infty)\longrightarrow
{\Cal{P}}_L(\infty)$ is {\it countably constructible} if it is of the form
$F^\infty_{\omega_{1}}$  where $F_{\omega_{1}}$ is a constructible function
from ${\Cal{P}}_L(\omega_1)$ to itself.

\proclaim{Theorem 14} Assume that $O^\#$ exists and $F:\
{\Cal{P}}_L(\infty)\longrightarrow {\Cal{P}}_L(\infty)$  is countably
constructible. Then there exists $R<_L O^\#$ such that $F(A)$ is definable
over $\langle L[R],A\rangle$  for all $A\in {\Cal{P}}_L(\infty).$ 
\endproclaim

\demo{Proof} Apply Theorem 13 to $F_{\omega_{1}}$ where
$F=F^\infty_{\omega_{1}}.$  The resulting real $R$  satisfies the conclusion
of the present Theorem. \hfill{$\dashv$}
\enddemo

\vskip5pt

\flushpar
{\bf Remarks.} \ (a) \ The definitions of $F(A)$ over $\langle
L_{\omega_{1}}(R),A\rangle,$ $\langle L[R],A\rangle$  
in Theorems 13, 14 respectively are
independent of $A.$ 

(b) \ If $F:\ {\Cal{P}}_L(\omega_1)\longrightarrow {\Cal{P}}_L(\omega_1)$  is
constructible then there exists a set-generic extension of $L$  in which there
is a real $R$  obeying the conclusion of Theorem 13. However we cannot expect
there to be such a real in $L[O^\#],$  or even compatible with the existence
of $O^\#.$  The key feature of our forcing ${\Cal{P}}$ is that not only can it
be used to produce a real $R$  obeying the conclusion of Theorem 12 but such a
real can be found in $L[O^\#].$  If one is willing to entirely ignore
compatibility with $O^\#$ then there are forcings far simpler than ours which
achieve the effect of Theorem 14 for any $F:$ classes  $\longrightarrow$
classes, over any model of G\"odel-Bernays class theory.

\enddocument